\documentclass[reqno,12pt]{amsart}
\usepackage{graphicx}
\usepackage{epsfig}
\usepackage{amssymb,amsthm,amsmath,amsfonts,amstext,mathrsfs,fancybox}
\usepackage{enumerate}
\usepackage{asymptote}
\usepackage{latexsym}
\usepackage{epsfig}
\usepackage{url}
\usepackage{fancyhdr}
\usepackage{multirow}
\usepackage{wrapfig}

\newtheorem*{thmm}{Theorem}

\newtheorem{prop}{Proposition}

\newtheorem*{defin}{Definition}

\newtheorem*{conj}{Conjecture}

\makeatletter
\def\dual#1{\expandafter\dual@aux#1\@nil}
\def\dual@aux#1/#2\@nil{\begin{tabular}{@{}l@{}}#1\\#2\end{tabular}}
\makeatother

\input xy
\xyoption{all}

\title[Full Grid Lattice Polygons]
{Full Grid Lattice Polygons\\ 
with Maximal Sum of Squares of Edge-Lengths}
\author[O. M. Ali\v{s}auskas, G. Alkauskas, V. Di\v{c}i\={u}nas]{Oliver Mantas Ali\v{s}auskas, Giedrius Alkauskas, Valdas Di\v{c}i\={u}nas}
\address{Vilnius University, Department of Mathematics and Informatics, Naugarduko 24, LT-03225 Vilnius, Lithuania}
\email{giedrius.alkauskas@mif.vu.lt; valdas.diciunas@mif.vu.lt}

\begin{document}
	\begin{abstract}
	Consider a subset $[1,2,\ldots,n]\times[1,2,\ldots,n]$ of the plane integer lattice. Take any non self-intersecting $n^2$-gon built on it (straight angles are allowed). The square of a side length is a positive integer. It is thus natural to ask how large the sum of square lengths of such an $n^2$-gon can be. This maximal value is a new integer sequence, labeled by A358212 in OEIS. In this note we give the lower bound and conjecture that this in fact is the correct answer. We further investigate proper $n^2$-gons (straight angles are not allowed) and present analogous results. Both sequences (conjecturally) have a different growth size.   
	\end{abstract}

\maketitle

\section{Full grid lattice polygons}
\subsection{Embedded graph} 

Given a positive integer $n\in\mathbb{N}$. Consider an integer lattice $\mathbb{Z}^{2}$ and its subset $\mathcal{A}_{n}=\{(x,y):(x,y)\in\mathbb{Z}^{2},\, 1\leq x,y\leq n\}$. Let $a_{i}=(x_{i},y_{i})$, $1\leq i\leq n^2$, be any ordering of this set. We can interpret it as a complete embedded graph $K_{n^2}$, the segment $e_{i,j}=[(x_{i},y_{i}),(x_{j},y_{j})]$, $i\neq j$, being an edge connecting two distinct vertices. The most natural way to define a weight of each particular edge is to count amount of unit moves needed to reach one vertex from the other. This definition does not dependent on any properties of the embedding: 
\begin{eqnarray*}
s_{0}(i,j)=|x_{i}-x_{j}|+|y_{i}-y_{j}|.
\end{eqnarray*}
However, in our case another obvious alternative is to consider weight function attached to an Euclidean geometry of the plane. We thus arrive to the definition       
\begin{eqnarray*}
	s(i,j)=(x_{i}-x_{j})^2+(y_{i}-y_{j})^2.
\end{eqnarray*}
Since now the geometry is involved, it makes sense to pose the following questions (m.w. stands for ``maximal weight", and n-i. stands for ``non-intersecting"):
\begin{itemize}
\item[1)] find m.w. Hamiltonian cycle with n-i. edges ($n^2$-gons supported on $\mathcal{A}_{n}$);
\item[2)]the same as above, additionally requiring non-congruence of any two consecutive edges;
\item[3)] m.w. n-i. spanning tree; 
\item[4)] m.w. n-i. triangulation.
\end{itemize}
In this paper we deal with questions 1) and 2). 
\subsection{Two sequences}

\begin{figure}
	\includegraphics[scale=0.24]{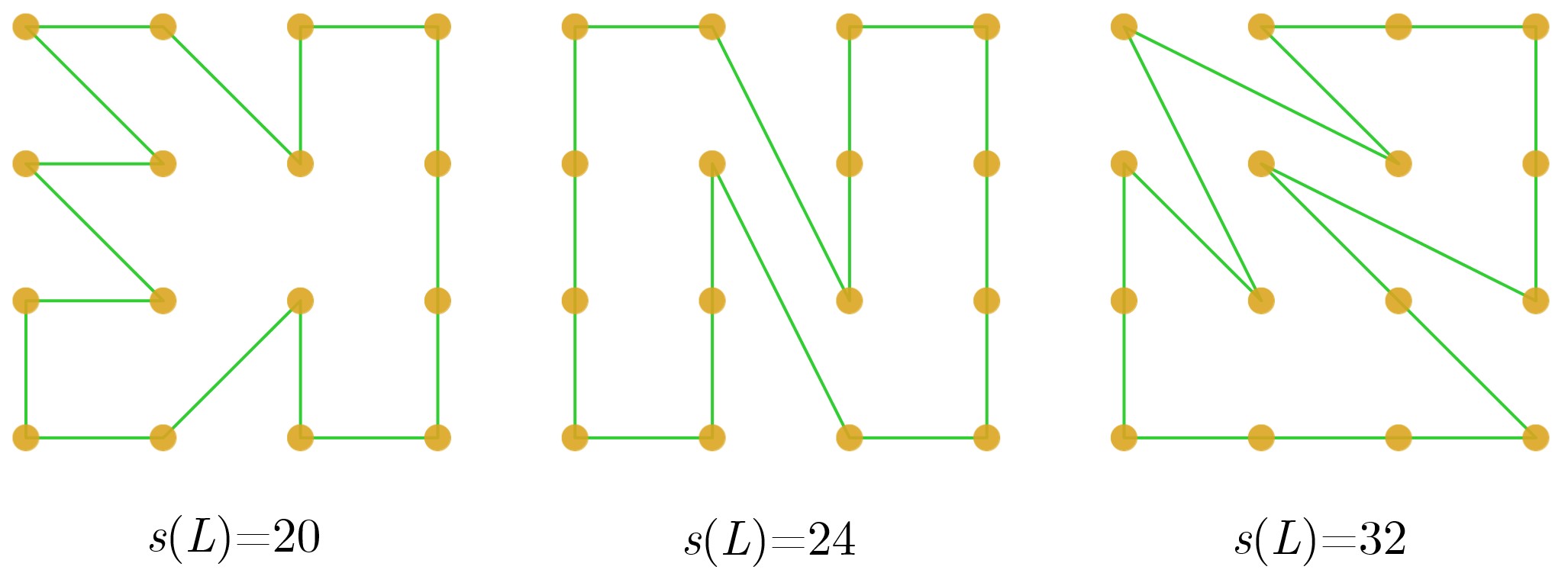}
	\caption{Three different closed broken lines for $n=4$}
	\label{pvz}
\end{figure}

And so, consider a closed non self intersecting broken line $L=(a_{1},a_{2},\ldots,a_{n^{2}})$ supported on $\mathcal{A}_{n}$.  We put $a_{n^2+1}=a_{1}$. Let

\begin{eqnarray}
	s(L)=\sum\limits_{i=1}^{n^{2}}|a_{i+1}-a_{i}|^{2}.\label{suma}
\end{eqnarray}
Here a point $a_{i}=(x_{i},y_{i})$ is treated as a $2$-dimensional vector. For $n=4$, three different examples are presented in Figure \ref{pvz}. 
\begin{defin}[Two main sequences]We define a sequence $a(n)$ by setting $a(n)$ to be the maximal possible value among $s(L)$ for all such $L$. Equally, let $a_{0}(n)$ be the maximum value of $s(L)$ among proper $n^2$-gons (straight angles are not allowed).
\end{defin}

\begin{table}
	\begin{tabular}{|c||r | r | r | r| r| r| r| r| r| r| r| r| r|}
		\hline
		\multicolumn{14}{|c|}{\textbf{Sequence $a(n)$}}\\
		\hline
		$n$&$2$&$3$&$4$&$5$&$6$&$7$&$8$&$9$ &$10$ &$11$ &$12$& $13$ & $14$\\
		\hline
		$a(n)$&$4$&$10$&$36$&$98$&$232$&$\mathbf{462}$&$\mathbf{842}$&$\mathbf{1424}$ &$\mathbf{2242}$ & $\mathbf{3378}$ & $\mathbf{4920}$ & $\mathbf{6906}$ & $\mathbf{9436}$\\
		Fjords&-&-&-&$(1,1)$&$(1,1)$&$(1,2)$&$(2,2)$&$(2,2)$ &$(2,2)$ & $(2,3)$ or $(3,3)$ & $(3,3)$ & $(3,3)$ & $(3,4)$\\
		\hline
	\hline
	\end{tabular}
	\caption{Initial terms, true and conjectural (in bold)}
	\label{table1}
\end{table}

\begin{table}
	\begin{tabular}{|c||r | r | r | r| r|r|}
		\hline
		\multicolumn{7}{|c|}{\textbf{Sequence $a(n)$}}\\
		\hline
		$n$&$15$&$16$&$17$&$18$&$19$&$20$\\
		\hline
		$a(n)$&$\mathbf{12638}$&$\mathbf{16560}$&$\mathbf{21318}$&$\mathbf{27066}$&$\mathbf{33884}$&$\mathbf{41884}$\\
		Fjords&$(4,4)$&$(4,4)$&$(4,5)$&$(5,5)$&$(5,5)$&$(5,6)$\\
		\hline
	\end{tabular}
	\caption{Next conjectural terms}
	\label{table2}
\end{table}

\begin{table}
	\begin{tabular}{|c||r | r | r | r| r| r| r| r| r| r| r| r| r|}
		\hline
		\multicolumn{14}{|c|}{\textbf{Sequence $a_{0}(n)$}}\\
		\hline
		$n$&$2$&$3$&$4$&$5$&$6$&$7$&$8$&$9$ &$10$ &$11$ &$12$& $13$ & $14$\\
	
		\hline
		$a_{0}(n)$&$4$&$0$&$20$&$0$&$142$&$\mathbf{346}$&$\mathbf{656}$&$\mathbf{1180}$ &$\mathbf{1808}$ & $\mathbf{2810}$ & $\mathbf{3552}$ & $?$ & $?$\\
		\hline
		\hline 
	\end{tabular}
	\caption{Initial terms and bounds}
	\label{table3}
\end{table}
\begin{figure}
	\includegraphics[scale=0.29]{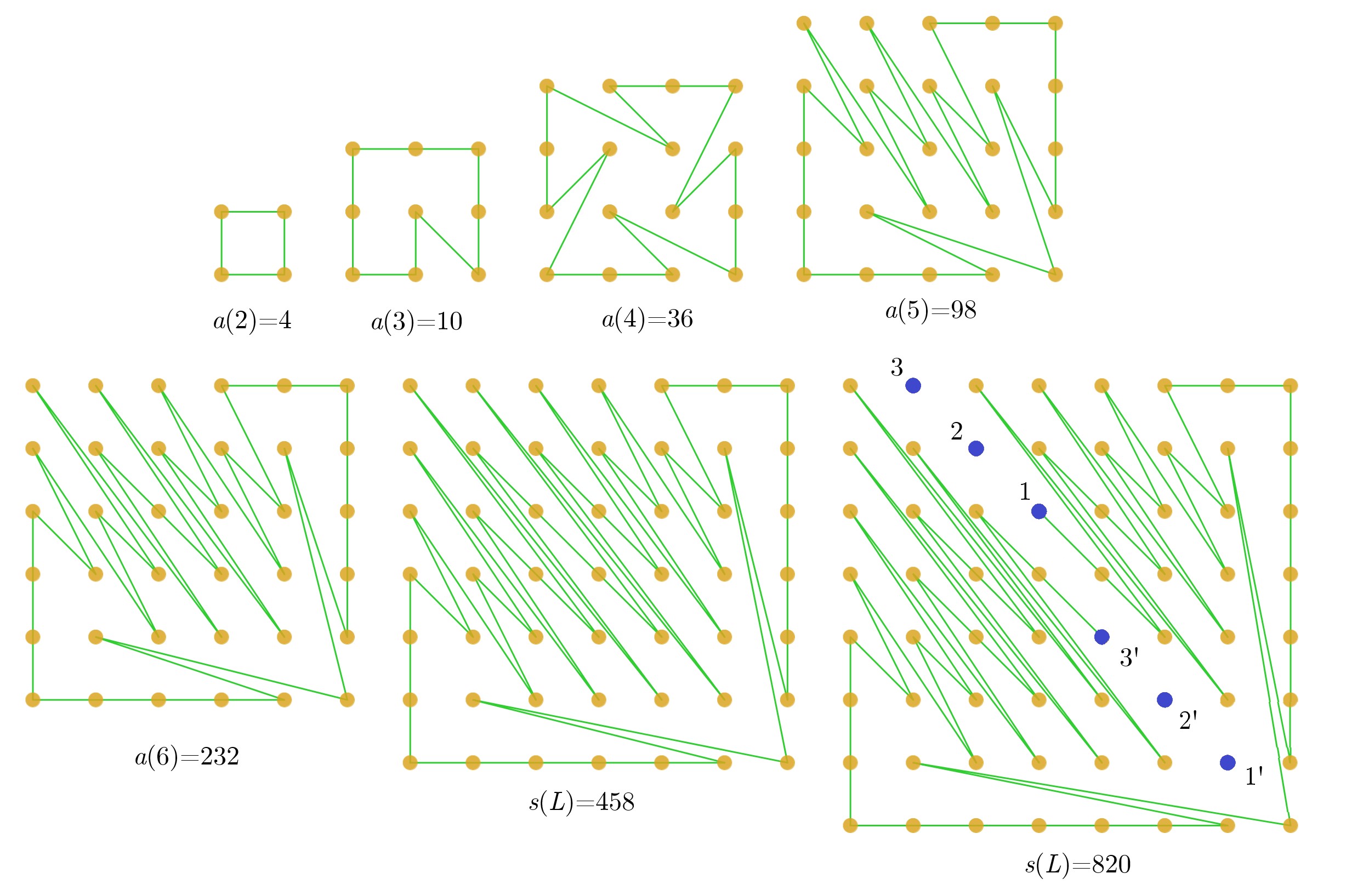}
	\caption{Extremal examples for the sequence $a(n)$ in cases $n=2$ through $6$ and the first approximation in cases $n=7,8$. It took several weeks of computer calculations to establish the exact value of $a(6)$. (In the last example the broken line goes $1$-$1'$-$2$-$2'$-$3$-$3'$).}
	\label{pirm-ast}
\end{figure}
\begin{figure}
	\includegraphics[scale=0.20]{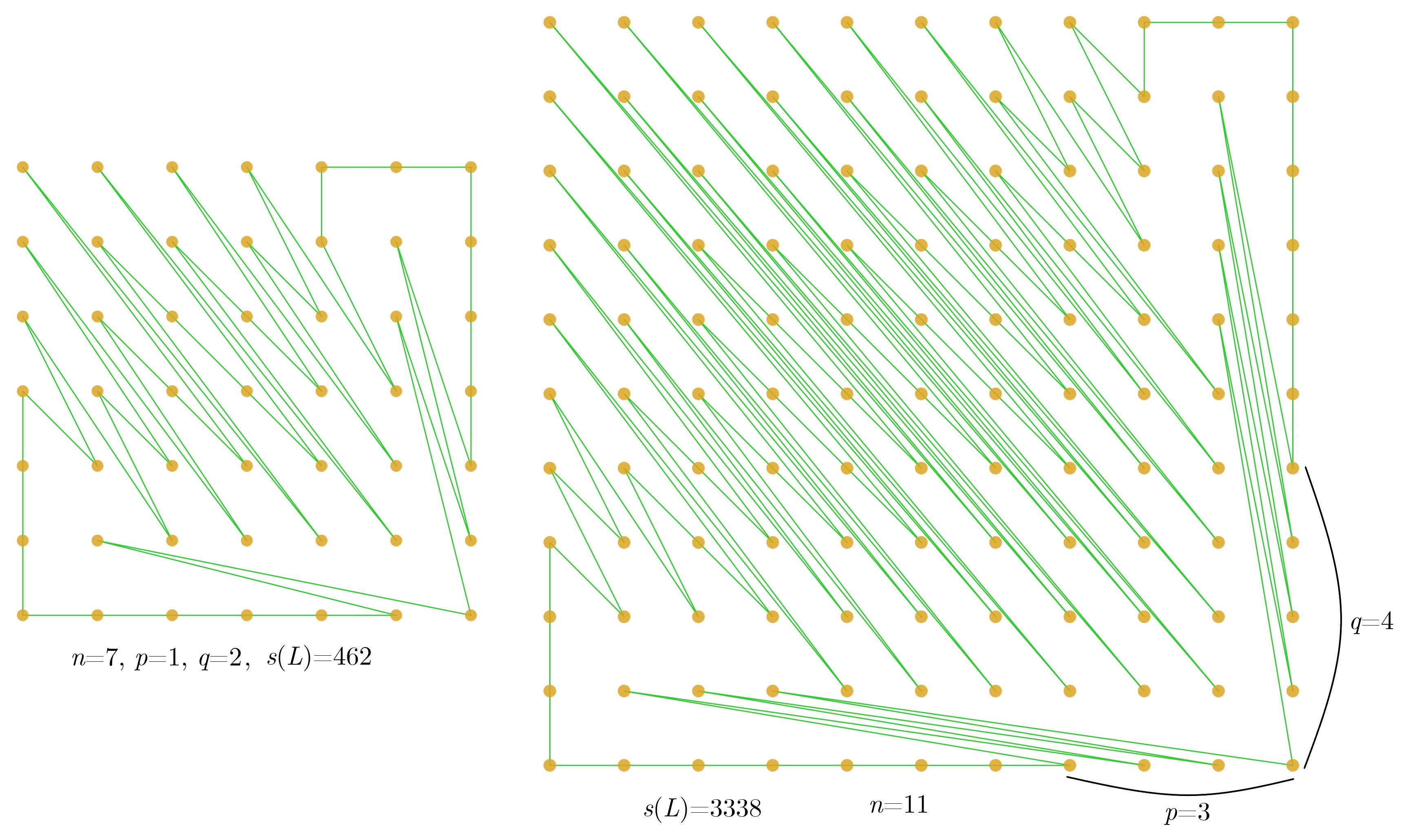}
	\caption{The optimal example for $n=7$ (left) and illustration of the notion of a ``fjord" ($n=11$, right).}
	\label{sept-vien}
\end{figure}
\begin{figure}
	\includegraphics[scale=0.16]{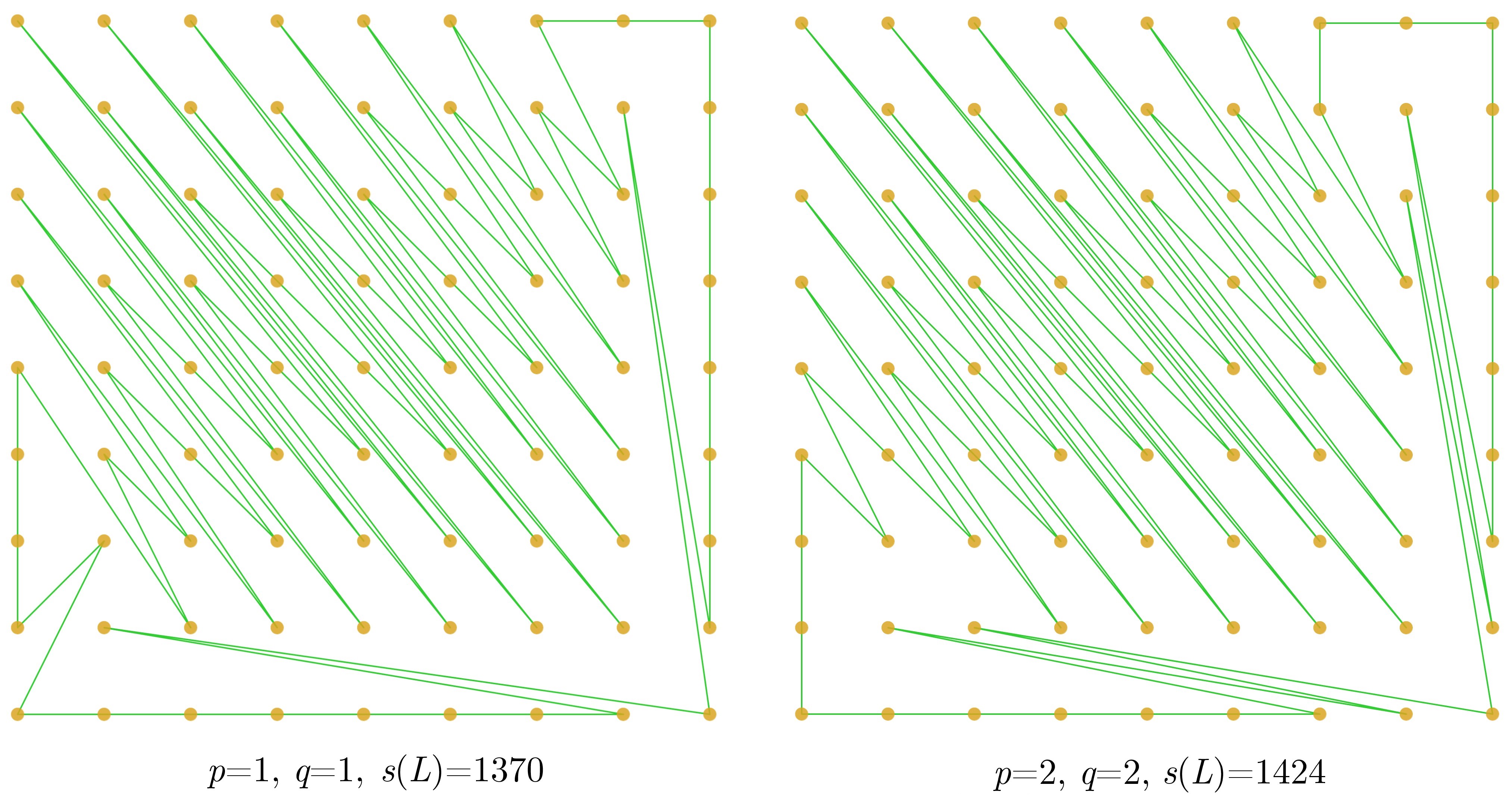}
	\caption{The first (left) and the second (right) approximations for $n=9$.}
	\label{devyni-comp}
\end{figure} 

\begin{figure}
	\includegraphics[scale=0.185]{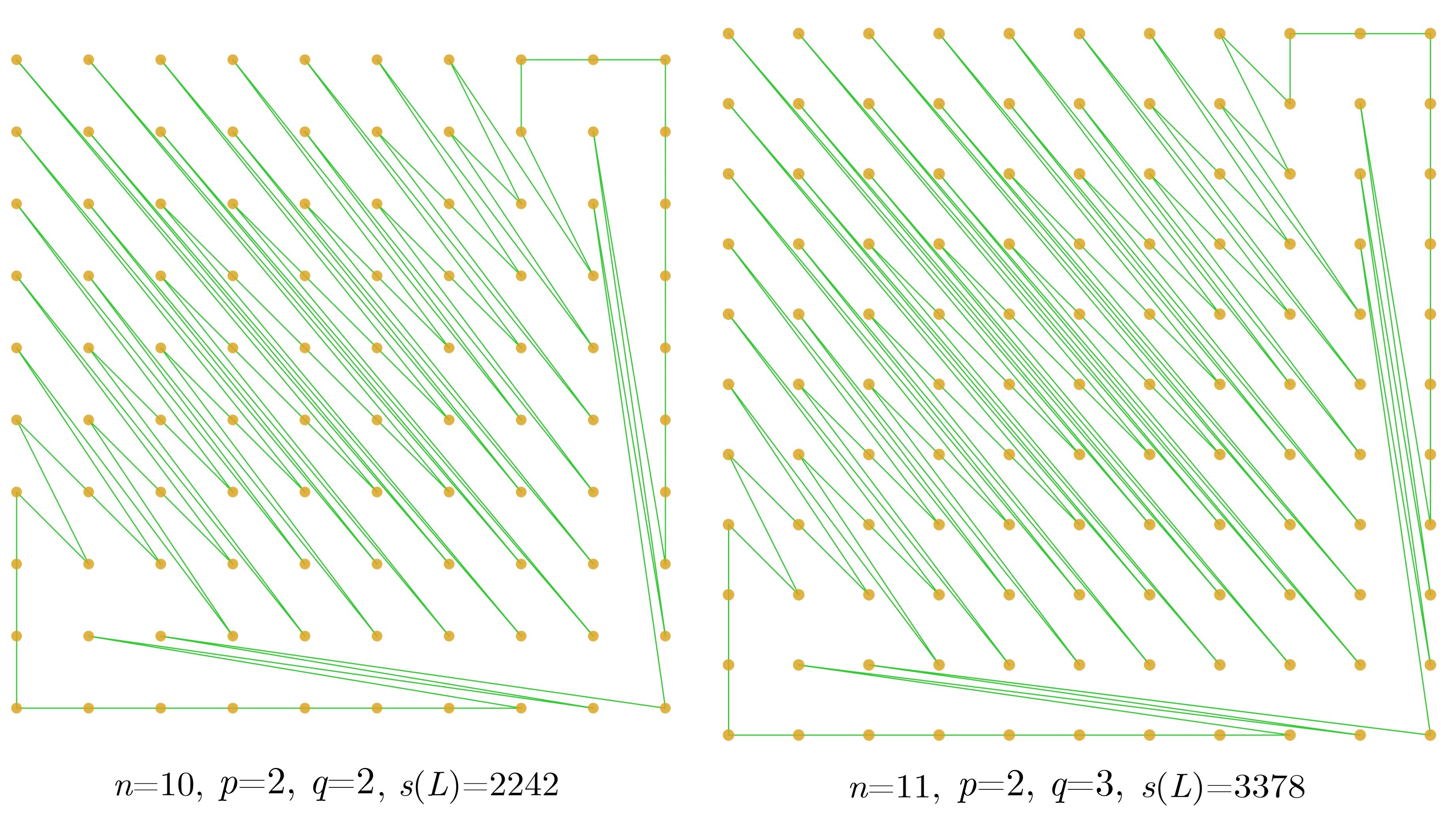}
	\caption{The second approximation for $n=10, 11$.}
	\label{desimt-vien}
\end{figure} 
The first sequence is labelled A358212 in OEIS \cite{oeis}. The initial exact and conjectural terms of these sequences are given in Tables \ref{table1},  \ref{table2} and \ref{table3}. Note that proper $n^2$-gons exist only for $n=4$ and $n\geq 6$ (see \cite{cgg}), so the sequence $a_{0}$ has an offset $n=6$.\\

The main result of this paper can now be immediately formulated. 
\begin{thmm}
	Let us define a sequence $\alpha(n)$, $n\in\mathbb{N}$, by setting	
		\begin{eqnarray*}
		\def\arraystretch{2.0}
		\alpha(n)=\frac{8}{27}n^4-\frac{464}{729}n^3+\left\{\begin{array}{l@{\qquad}r}
			-\,\,\,\,\,\frac{4}{3}n^2+\,\,\,\,6n-\,\,\,\,\,\,2,\, n\equiv 0;\quad
			-\frac{292}{243}n^2+\frac{922}{243}n-\frac{1642}{729},\, n\equiv 1;\\
			-\frac{308}{243}n^2+\frac{946}{243}n-\frac{1724}{729},\, n\equiv 2;\quad
			-\,\,\,\,\frac{4}{3}n^2+\frac{158}{27}n-\,\,\,\,
			\frac{64}{27},\, n\equiv 3;\\
			-\frac{292}{243}n^2+\frac{886}{243}n-\frac{1804}{729},\,
			 n\equiv 4;\quad
			-\frac{308}{243}n^2+\frac{874}{243}n-\frac{1400}{729},\, n\equiv 5;\\
			-\,\,\,\,\frac{4}{3}n^2+\frac{158}{27}n-\,\,\,\,\frac{44}{27},\, n\equiv 6;\quad
			-\frac{292}{243}n^2+\frac{886}{243}n-\frac{1264}{729},\, n\equiv 7;\\
			-\frac{308}{243}n^2+\frac{946}{243}n-\frac{1292}{729},\, n\equiv 8.

		\end{array}\right.
	\end{eqnarray*}
here $n\equiv i$ is an abbreviation of a clause ``if $n\equiv i\,\mathrm{ (mod}\,9)$". \\

Then for $n\geq 9$, one has a bound $a(n)\geq \alpha(n)$.
	\end{thmm}
This form of the result is most suitable if one is interested in asymptotics. For example, if all $9$ polynomials are ordered lexicographically, then the subcase $n\equiv i\,\mathrm{ (mod}\,9)$ has a slight growth edge over the rest. However, the structure of a broken line $L$ is much more transparent if the result is stated as in Proposition \ref{antras}.
\subsection{Development}Our work on this problem consists of three ingredients:
\begin{itemize}
	\item[$\star$]Computer search. Though this part is invisible in our text, computations lasted for several months; gradual refinements of the algorithms involved were implemented.
	\item[$\star$]A theoretical part (the rest body of the paper).
	\item[$\star$]A web application (available at \cite{app}) for visualization, data import-export (needed for MAPLE calculations), and for double-checking all the theoretical results.
	\end{itemize} 
The example for $n=6$ was verified via a complete search. This suggested the first approximation to our problem (Subsection \ref{first-apprx}, Figure \ref{pirm-ast}). However, already in the case $n=7$ computer found that, building on the first approximation, an can be made (Figure \ref{sept-vien}, left). This is achieved with the help of a construction which we label as ``fjords". Figure \ref{sept-vien} (right) provides a visual explanation of this notion. In Subsection \ref{second-apprx} we calculate the optimal number of fjords needed to maximize the sum $s(L)$.\\

At this point, no immediate local improvement to the construction can be contrived. This allows us to cautiously pose the following Conjecture. 
\begin{conj} For $n\geq 9$, $a(n)=\alpha(n)$.
\end{conj}
  
\section{Bounds}

\subsection{Lower bound}
\begin{figure}
	\includegraphics[scale=0.17]{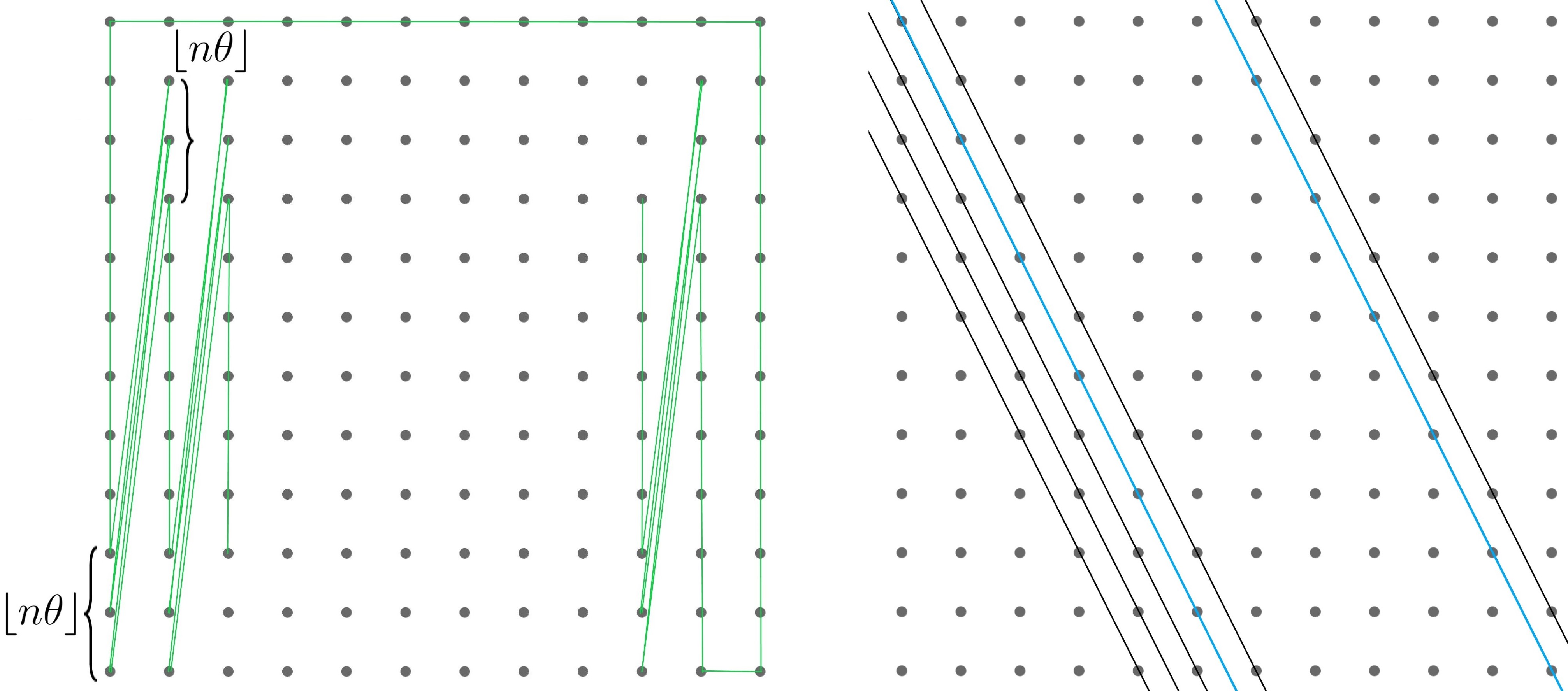}
	\caption{The first construction (left); dissection of $\mathcal{A}_{n}$ into collinear sets (right)}
	\label{lower}
\end{figure}

Let $\theta\in\mathbb{R}$, $0<\theta<\frac{1}{2}$. To establish the simple lower bound, consider Figure \ref{lower} (left). The construction of a polygon is visually apparent. 
 Direct calculation gives
\begin{eqnarray}
	\begin{split}
	s(L_{n})&=&\Big{(}\big{(}n-\lfloor n\theta\rfloor\big{)}^2+1\Big{)}\cdot(\lfloor n\theta\rfloor-1)\cdot (n-2)\\
	&+&\Big{(}\big{(}n-\lfloor n\theta\rfloor-1\big{)}^2+1\Big{)}\cdot \lfloor n\theta\rfloor \cdot (n-2)\\
	&+&(n-2\lfloor n\theta\rfloor)\cdot(n-3)+4n-2\lfloor n\theta\rfloor-2\\
	&=&2n^4\cdot (1-\theta)^2\theta+O(n^3).\label{pirmun}
	\end{split}
	\end{eqnarray}
The function $\theta(1-\theta)^2$ achieves a local maximum $\frac{4}{27}$ at $\theta_{0}=\frac{1}{3}$. \emph{A posteriori}, this simple fact is solely responsible for the leading term in the asymptotics. Thus, we obtain
	\begin{eqnarray*}
		\liminf\limits_{n\rightarrow\infty}\frac{a(n)}{n^4}\geq\frac{8}{27}.
	\end{eqnarray*}
A more careful calculation in (\ref{pirmun}) for $\theta=\theta_{0}$ shows that
\begin{eqnarray*}
	s(L_{n})=\frac{8}{27}n^4-\frac{40}{27}n^3+O(n^2).
	\end{eqnarray*}

The presented construction is based on the idea that $\mathcal{A}_{n}$ is fully covered by the union of $n$ lines $\mathcal{L}_{i}=\{(x,y)\in\mathbb{R}^2:x=i\}$, $i=1,2,\ldots, n$. Each of them contains $n$ grid points. In general, the set $\mathcal{A}_{n}$ can be dissected into parallel collinear finite sets in many ways. For example,  Figure \ref{lower} (right) gives such a partition for the direction $\overrightarrow{p}=(-1,2)$. The set $\mathcal{A}_{n}$ is thus split into $4n-3$ collinear sets.\\

 We will now focus on the case $\overrightarrow{p}=(-1,1)$. The task is to generalize the pattern shown in Figure \ref{pirm-ast} for any $n\geq 9$.
 \section{The sharper lower bound}
\subsection{The first approximation}
\label{first-apprx}

As a first step in proving Theorem, we will demonstrate the validity of the following result.
\begin{prop} Define the sequence $\beta(n)$ by
	\begin{eqnarray*}
		\def\arraystretch{2.0}
		\beta(n)=\frac{8}{27}\,n^4-\frac{32}{27}\,n^3+\frac{16}{3}\,n^2-18n+20+\left\{\begin{array}{r@{\qquad}r}
			0,\text{ if } n\equiv 0\,\mathrm{ (mod}\,3),\\
			\displaystyle{-\frac{8n+4}{27}},\text{ if } n\equiv 1\,\mathrm{ (mod}\,3),\\
			\displaystyle{-\frac{16n-16}{27}},\text{ if } n\equiv 2\,\mathrm{ (mod}\,3).
		\end{array}\right.
	\label{pagr}
	\end{eqnarray*}
	Then $a(n)\geq \beta(n)$.
	\label{prop1}
\end{prop}
	

	\begin{figure}
		\includegraphics[scale=0.2]{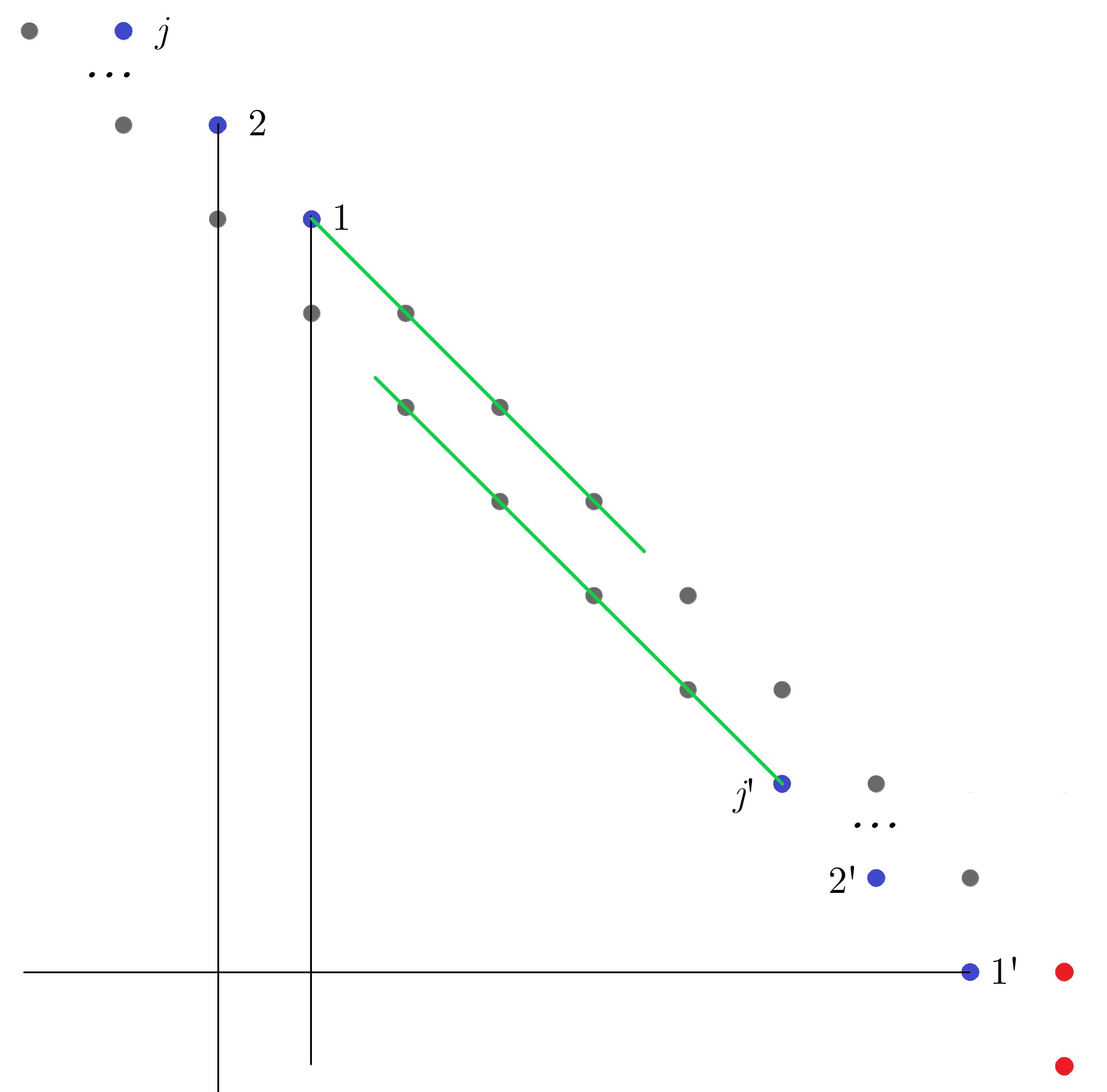}
		\caption{Main construction}
		\label{recu}
	\end{figure}
To prove this, let us first take a non-principal diagonal with $k$ points ($k\leq n-1$) in the upper triangle of $\mathcal{A}_{n}$ (the one above the main diagonal $x=y$). Mark its first $j$ upper points. In the same manner, let us mark the last but one $j$ points of the $(k+1)$th diagonal. The whole construction is presented in Figure \ref{recu}. Consider now the broken line $1$-$1'$-$2$-$2'$-$\cdots$-$j$-$j'$ as a part of our $n^2$-gon. The task is to maximize the contribution from this broken line to the sum (\ref{suma}).  The gain function is given by
\begin{eqnarray*}
F&=&\Big{[}(k-j-1)^2+(k-j)^2\Big{]}j+\Big{[}(k-j)^2+(k-j+1)^2\Big{]}(j-1)-4(j-1)\\
	&=&4j^3+(-8k-2)j^2+(4k^2+4k)j-2k^2-2k+3.
	\end{eqnarray*}

The first two summands are actually sum of squares of side lengths ($j$ shorter and $(j-1)$ longer hypotenuses of right triangles), while the term $-4(j-1)$ indicates that each additional blue point shortens both green lines. Now, $\frac{d}{d j}\,F(j)=4(3j-k-1)(j-k)$.
Thus, the maximum is achieved for $j=\frac{k+1}{3}$. If this is not an integer, the value of $F(k,j)$ at both points $j=\big{\lfloor}\frac{k}{3}\big{\rfloor}$ and $j=\big{\lfloor}\frac{k}{3}\big{\rfloor}+1$ must be compared. Simple calculations yield the optimal value
	$j=\big{\lceil}\frac{k}{3}\big{\rceil}$. And so, the maximal gain is given by
\begin{eqnarray*}
	\def\arraystretch{2.0}
	U(k)=\max\limits_{j\in\mathbb{N},\, 1\leq j<\frac{k}{2}}F(j,k)=\left\{\begin{array}{c@{\quad}l}
	\displaystyle{F\Big{(}k,\frac{k}{3}\Big{)}}=16s^3-8s^2-6s+3,&\text{ if }k=3s;\\
		\displaystyle{F\Big{(}k,\frac{k+2}{3}\Big{)}}=16s^3+8s^2-6s+1,&\text{ if } k=3s+1;\\
		\displaystyle{F\Big{(}k,\frac{k+1}{3}\Big{)}}=16s^3+24s^2+6s+1,
		&\text{ if } k=3s+2.
	\end{array}\right.
\end{eqnarray*}

The lower triangle case is dealt with similarly. Let $\beta(n)=s(L)$ be the corresponding sum maximized over all diagonals. All these calculations yield a recurrence
\begin{eqnarray}
	\beta(n)=\beta(n-1)+2U(n-1)+12n-32,\quad n\geq 6. \label{perein}
\end{eqnarray}
\emph{A priori}, it is clear that the true recurrence for $\beta$ is \emph{approximately} of this form. The best way to verify that it is \emph{exactly} of this form is not to give a tedious diagonal-by-diagonal summation, but to play around the files \cite{files}, uploaded into the application \cite{app}. \\

The recurrence (\ref{perein}), on its turn, gives the exact formula in Proposition \ref{prop1}.

\subsection{The second approximation. ``Fjords"}
\label{second-apprx}
 As we will see now, the first approximation is not sharp for $n\geq 7$. This fact was discovered with the help of computer. The phenomenon behind this improvement can be seen in the Figure \ref{devyni-comp}. In words: the construction in the first approximation can be locally altered to yield a better result. Our task to to find an optimal alteration.  
\begin{prop} 
	\label{antras}
Let $n\geq 9$, $1\leq p,q\leq n-3$. The described construction with $p$ fjords at the bottom and $q$ fjords to the right gives the sum 
\begin{eqnarray*}
s(n;p,q)&=&\beta(n)-4n^2+20n-18\\
&+&p(n-p-2)^2+p(n-p-1)^2-U(p+2)-p\\
&+&q(n-q-2)^2+q(n-q-1)^2-U(q+1)-q.
	\end{eqnarray*}
The maximum is achieved when $p=\lfloor\frac{n-3}{3}\rfloor$ and, respectively,  $q=\lfloor\frac{n-2}{3}\rfloor$.
\end{prop}
Yet again, we can provide a tedious step-by-step derivation of this result. As noted, alternatively one can play around a set of files \cite{files} to convince oneself that the statement holds true. 
\section{The second sequence}
\label{second}
We cannot now use $\theta_{0}=\frac{1}{3}$ in the construction, since it gives many straight angles. Yet we still can use $\theta=\frac{1}{2}$. This avoids the problem of degeneration of $n^2$-gon, except for the border of the grid. Yet it is obvious that the following holds:
	\begin{eqnarray*}
	\liminf\limits_{n\rightarrow\infty}\frac{a_{0}(n)}{n^4}\geq\frac{1}{4}.
\end{eqnarray*}

\begin{figure}
	\includegraphics[scale=0.43]{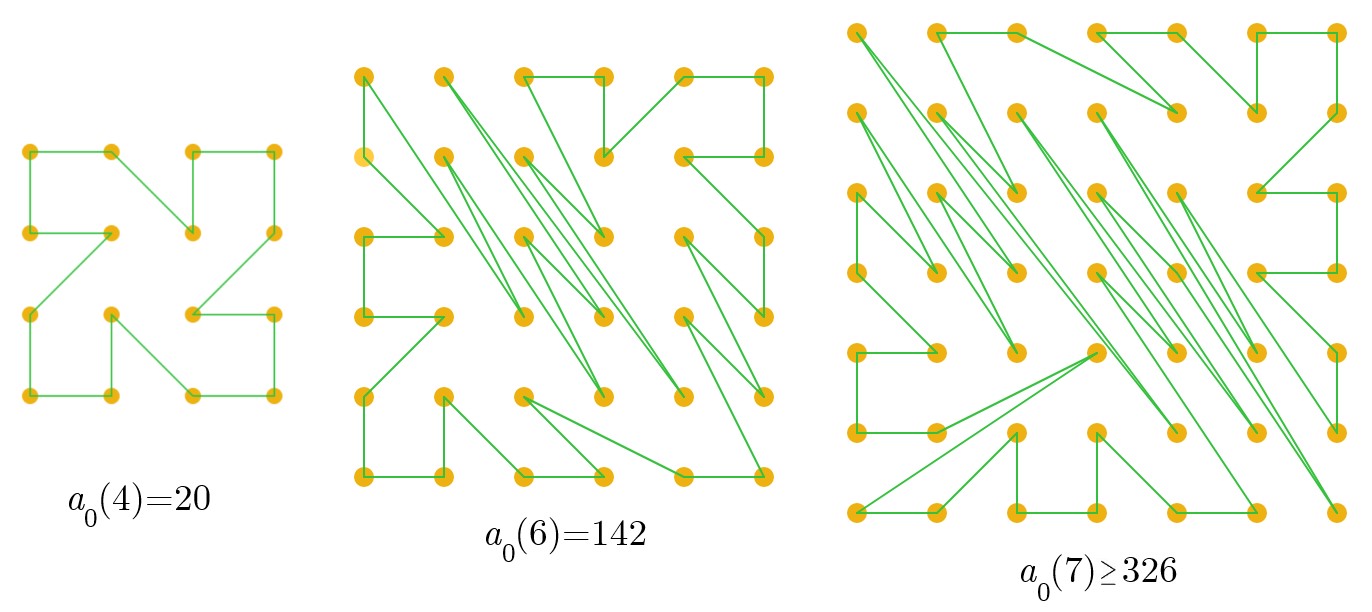}
	\caption{The sequence $a_{0}$. Extremal example for $n=4,6$ and (most probably) for $n=7$.}
\end{figure}
\begin{figure}
	\includegraphics[scale=0.15]{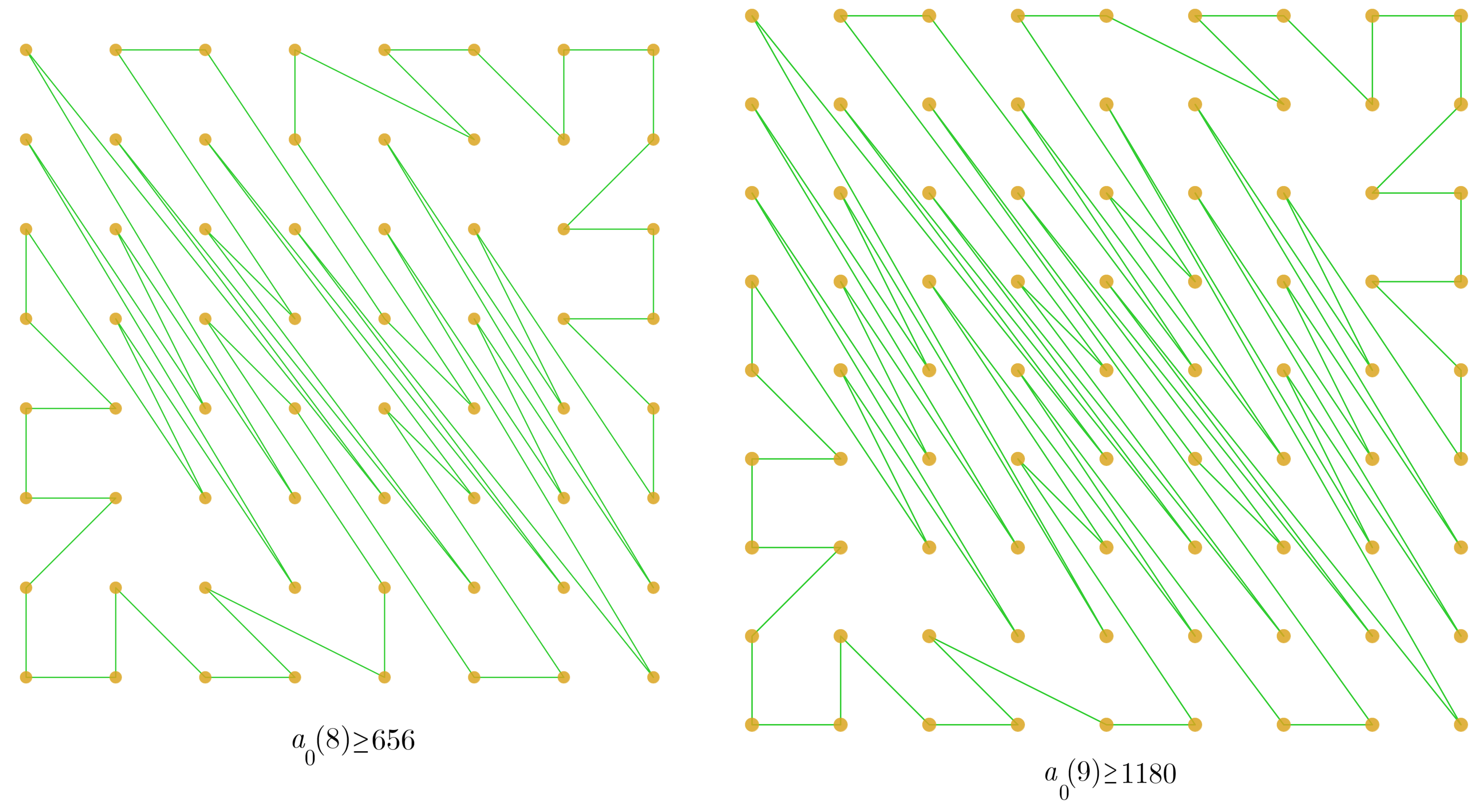}
	\caption{The sequence $a_{0}$. Best (so far) examples for $n=8,9$.}
\end{figure}
\begin{figure}
	\includegraphics[scale=0.17]{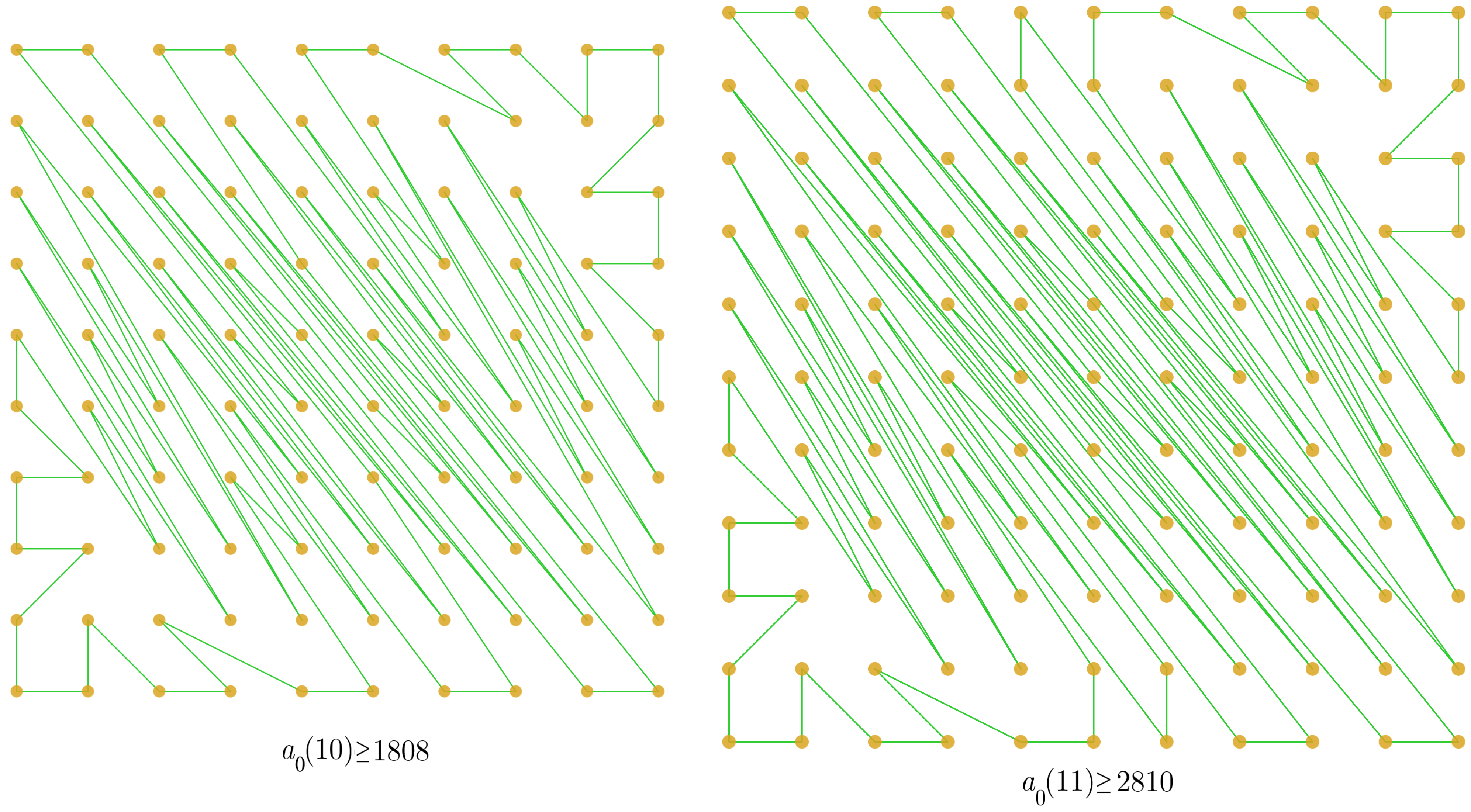}
	\caption{The sequence $a_{0}$. Best (so far) examples for $n=10,11$.}
\end{figure}


\subsection{The upper bound}
It is fairly easy to show that
\begin{eqnarray*}
	\limsup\limits_{n\rightarrow\infty}\frac{a(n)}{n^4}\leq \frac{2}{3}.
\end{eqnarray*}
Though this exceeds the conjectural bound more than twice, we unfortunatelly have not yet managed to improve on this. 

\end{document}